\input amstex

\documentstyle{amsppt}
  \magnification=1100
  \hsize=6.2truein
  \vsize=9.0truein
  \hoffset 0.1truein
  \parindent=2em

\NoBlackBoxes


\font\eusm=eusm10                   


\font\eusms=eusm7                       

\font\eusmss=eusm5                      


\newcount\theTime
\newcount\theHour
\newcount\theMinute
\newcount\theMinuteTens
\newcount\theScratch
\theTime=\number\time
\theHour=\theTime
\divide\theHour by 60
\theScratch=\theHour
\multiply\theScratch by 60
\theMinute=\theTime
\advance\theMinute by -\theScratch
\theMinuteTens=\theMinute
\divide\theMinuteTens by 10
\theScratch=\theMinuteTens
\multiply\theScratch by 10
\advance\theMinute by -\theScratch

\def\today{{\number\day\space
 \ifcase\month\or
  January\or February\or March\or April\or May\or June\or
  July\or August\or September\or October\or November\or December\fi
 \space\number\year}}

\define\Ad{\text{\rm Ad}}

\define\Afr{{\frak A}}

\define\biggnm#1{
  \bigg|\bigg|#1\bigg|\bigg|}

\define\bignm#1{
  \big|\big|#1\big|\big|}

\define\Bo{{B\oup}}

\define\Cpx{\bold C}

\define\dif{\text{\it d}}

\define\eqdef{{\;\overset\text{def}\to=\;}}

\define\Eto#1{E_{(\to{#1})}}

\define\fpamalg#1{{\dsize\;\operatornamewithlimits*_{#1}\;}}

\define\fpiamalg#1{{\tsize\;({*_{#1}})_{\raise-.5ex\hbox{$\ssize\iota\in I$}}}}

\define\freeprod#1#2{\mathchoice
     {\operatornamewithlimits{\ast}_{#1}^{#2}}
     {\raise.5ex\hbox{$\dsize\operatornamewithlimits{\ast}
      _{#1}^{#2}$}\,}
     {\text{oops!}}{\text{oops!}}}

\define\freeprodi{\mathchoice
     {\operatornamewithlimits{\ast}
      _{\iota\in I}}
     {\raise.5ex\hbox{$\dsize\operatornamewithlimits{\ast}
      _{\sssize\iota\in I}$}\,}
     {\text{oops!}}{\text{oops!}}}

\define\freeprodvni{\mathchoice
      {\operatornamewithlimits{\overline{\ast}}
       _{\iota\in I}}
      {\raise.5ex\hbox{$\dsize\operatornamewithlimits{\overline{\ast}}
       _{\sssize\iota\in I}$}\,}
      {\text{oops!}}{\text{oops!}}}

\define\Hil{{\mathchoice
     {\text{\eusm H}}
     {\text{\eusm H}}
     {\text{\eusms H}}
     {\text{\eusmss H}}}}


\define\Hilto#1{\Hil_{(\to{#1})}}

\define\Lambdao{{\Lambda\oup}}

\define\ld#1{{\hbox{..}(#1)\hbox{..}}}

\define\lrnm#1{\left\|#1\right\|}

\define\lspan{\text{\rm span}@,@,@,}

\define\nm#1{\|#1\|}

\define\Nats{\Naturals}

\define\Naturals{{\bold N}}

\define\otdts#1{\otimes_{#1}\cdots\otimes_{#1}}

\define\oup{^{\text{\rm o}}}

\define\owedge{{
     \operatorname{\raise.5ex\hbox{\text{$
     \ssize{\,\bigcirc\llap{$\ssize\wedge\,$}\,}$}}}}}

\define\owedgeo#1{{
     \underset{\raise.5ex\hbox
     {\text{$\ssize#1$}}}\to\owedge}}

\define\Pto#1{{P_{(\to{#1})}}}


\define\pup#1#2{{{\vphantom{#2}}^{#1}\!{#2}}\vphantom{#2}}

\define\QED{\newline
            \line{$\hfill$\qed}\enddemo}

\define\Reals{{\bold R}}

\define\restrict{\lower .3ex
     \hbox{\text{$|$}}}

\define\smd#1#2{\underset{#2}\to{#1}}

\define\smdb#1#2{\undersetbrace{#2}\to{#1}}

\define\smdbp#1#2#3{\overset{#3}\to
     {\smd{#1}{#2}}}

\define\smdbpb#1#2#3{\oversetbrace{#3}\to
     {\smdb{#1}{#2}}}

\define\smdp#1#2#3{\overset{#3}\to
     {\smd{#1}{#2}}}

\define\smdpb#1#2#3{\oversetbrace{#3}\to
     {\smd{#1}{#2}}}

\define\smp#1#2{\overset{#2}\to
     {#1}}

\define\Tcirc{\bold T}

\define\tocdots
  {\leaders\hbox to 1em{\hss.\hss}\hfill}


  \newcount\mycitestyle \mycitestyle=1 

  \newcount\bibno \bibno=0
  \def\newbib#1{\advance\bibno by 1 \edef#1{\number\bibno}}
  \ifnum\mycitestyle=1 \def\cite#1{{\rm[\bf #1\rm]}} \fi
  \def\scite#1#2{{\rm[\bf #1\rm, #2]}}


  \newcount\ignorsec \ignorsec=0
  \def\notasec{\ignorsec=1}

  \newcount\secno \secno=0
  \def\newsec#1{\procno=0 \subsecno=0 \ignorsec=0
    \advance\secno by 1 \edef#1{\number\secno}
    \edef\currentsec{\number\secno}}

  \newcount\subsecno
  \def\newsubsec#1{\procno=0 \advance\subsecno by 1
    \edef\currentsec{\number\secno.\number\subsecno}
     \edef#1{\currentsec}}

  \newcount\appendixno \appendixno=0
  \def\newappendix#1{\procno=0 \ignorsec=0 \advance\appendixno by 1
    \ifnum\appendixno=1 \edef\appendixalpha{\hbox{A}}
      \else \ifnum\appendixno=2 \edef\appendixalpha{\hbox{B}} \fi
      \else \ifnum\appendixno=3 \edef\appendixalpha{\hbox{C}} \fi
      \else \ifnum\appendixno=4 \edef\appendixalpha{\hbox{D}} \fi
      \else \ifnum\appendixno=5 \edef\appendixalpha{\hbox{E}} \fi
      \else \ifnum\appendixno=6 \edef\appendixalpha{\hbox{F}} \fi
    \fi
    \edef#1{\appendixalpha}
    \edef\currentsec{\appendixalpha}}

  \newcount\procno \procno=0
  \def\newproc#1{\advance\procno by 1
   \ifnum\ignorsec=0 \edef#1{\currentsec.\number\procno}
                     \edef\currentproc{\currentsec.\number\procno}
   \else \edef#1{\number\procno}
         \edef\currentproc{\number\procno}
   \fi}

  \newcount\subprocno \subprocno=0
  \def\newsubproc#1{\advance\subprocno by 1
   \ifnum\subprocno=1 \edef#1{\currentproc a} \fi
   \ifnum\subprocno=2 \edef#1{\currentproc b} \fi
   \ifnum\subprocno=3 \edef#1{\currentproc c} \fi
   \ifnum\subprocno=4 \edef#1{\currentproc d} \fi
   \ifnum\subprocno=5 \edef#1{\currentproc e} \fi
   \ifnum\subprocno=6 \edef#1{\currentproc f} \fi
   \ifnum\subprocno=7 \edef#1{\currentproc g} \fi
   \ifnum\subprocno=8 \edef#1{\currentproc h} \fi
   \ifnum\subprocno=9 \edef#1{\currentproc i} \fi
   \ifnum\subprocno>9 \edef#1{TOO MANY SUBPROCS} \fi
  }

  \newcount\tagno \tagno=0
  \def\newtag#1{\advance\tagno by 1 \edef#1{\number\tagno}}



\notasec
 \newtag{\BugA}
 \newtag{\AdenseBug}

\newbib{\Avitzour}
\newbib{\BercoviciVoiculescuZZSuperconv}
\newbib{\CuntzZZKthCertain}
\newbib{\DykemaZZFaithful}
\newbib{\DykemaZZPII}
\newbib{\DykemaHaagerupRordam}
\newbib{\DykemaRordamZZPI}
\newbib{\DykemaRordamZZProj}
\newbib{\DykemaRordamZZProjII}
\newbib{\KishimotoKumjian}
\newbib{\VoiculescuZZSymmetries}
\newbib{\VDNbook}

\topmatter
  \title Purely Infinite, Simple $C^*$-algebras Arising from Free
         Product Constructions, III
  \endtitle

  \author Marie Choda and Kenneth J\. Dykema \endauthor

  \date \today \enddate

  \rightheadtext{}

  \leftheadtext{}

  \address Dept.~of Mathematics,
           Osaka Kyoiku University,
           Asahigaoka, Kashiwara 582, Japan
  \endaddress

  \email marie\@cc.osaka-kyoiku.ac.jp 
  \endemail

  \address Dept.~of Mathematics,
           Texas A\&M University,
           College Station TX 77843--3368, USA
  \endaddress

  \email Ken.Dykema\@math.tamu.edu, {\it Internet:}
         http://www.math.tamu.edu/\~{\hskip0.1em}Ken.Dykema/
  \endemail

  \abstract
    In the reduced free product of C$^*$--algebras,
    $(A,\phi)=(A_1,\phi_1)*(A_2,\phi_2)$ with respect to faithful states
    $\phi_1$ and $\phi_2$, $A$ is  purely
    infinite and simple if $A_1$ is a reduced crossed product
    $B\rtimes_{\alpha,r}G$ for $G$ an infinite group, if $\phi_1$ is well
    behaved with respect to this crossed product decomposition, if $A_2\ne\Cpx$
    and if $\phi$ is not a trace.
  \endabstract

  \subjclass 46L05, 46L35 \endsubjclass

\endtopmatter

\document \TagsOnRight \baselineskip=18pt

  The reduced free product construction for C$^*$--algebras was invented
independently by Voiculescu~\cite{\VoiculescuZZSymmetries} and, in a more
limited sense, Avitzour~\cite{\Avitzour}.
(The term ``reduced'' is to distinguish this construction from the universal or
``full'' free product of C$^*$--algebras.)
It is a natural construction in Voiculescu's free probability theory,
(see~\cite{\VDNbook}).
Given unital C$^*$--algebras $A_\iota$ with states $\phi_\iota$ whose GNS
representations are faithful, ($\iota\in I$), the construction yields
$$ (A,\phi)=\freeprodi(A_\iota,\phi_\iota), $$
where $A$ is a unital C$^*$--algebra containing copies
$A_\iota\hookrightarrow A$ and generated by $\bigcup_{\iota\in I}A_\iota$, and
where $\phi$ is a state on $A$ with faithful GNS representation that restricts
to give $\phi_\iota$ on $A_\iota$ for every $\iota\in I$ and such that
$(A_\iota)_{\iota\in I}$ is free with respect to $\phi$.
Moreover, $\phi$ is a trace if and only if every $\phi_\iota$ is a trace;
by~\cite{\DykemaZZFaithful}, $\phi$ is faithful on $A$ if and only if
$\phi_\iota$ is faithful on $A_\iota$ for every $\iota\in I$.

It is a very interesting open question whether every simple, unital C$^*$--algebra
must either have a trace or be purely infinite.
Purely infinite C$^*$--algebras were defined by
J\. Cuntz~\cite{\CuntzZZKthCertain}.
A simple unital C$^*$--algebra $A$ is purely infinite if and only if for every
positive element $x\in A$ there is $y\in A$ with $y^*xy=1$.
An equivalent condition is that every hereditary C$^*$--subalgebra of $A$
contains an infinite projection.

Let
$$ (A,\phi)=(A_1,\phi_1)*(A_2,\phi_2) $$
be a reduced free product of C$^*$--algebras.
In~\cite{\DykemaRordamZZProj} it was shown that if $\phi_1$ or $\phi_2$ is
nontracial and if $A_1$ and $A_2$ are not too small in a specific sense, then
$A$ is properly infinite.
It is a plausible conjecture that whenever $A$ is simple and at least one of
$\phi_1$ and $\phi_2$ is not a trace, the C$^*$--algebra $A$ must be purely
infinite.
The first results in this direction were~\cite{\DykemaRordamZZPI}, where in a
certain class of examples when $\phi_1$ was assumed to be non--faithful, $A$
was shown to be purely infinite and simple.
In~\cite{\DykemaZZPII}, assuming $\phi_1$ and $\phi_2$ faithful, $A$ was shown
to be purely infinite and simple in the case when the centralizer of $\phi_1$
in $A_1$ contains a diffuse abelian subalgebra and when $A_2$ contains a
partial isometry that, loosely
speaking, scales $\phi_2$ by a constant $\lambda\ne1$.
In~\cite{\DykemaRordamZZProjII}, reduced free products of (countably)
infinitely many C$^*$--algebras that are not too small in a specific sense were
shown to be purely infinite.

In this note, we prove a theorem implying that $A$ is purely infinite and
simple under somewhat different conditions.
For example, if $A_1=C(\Tcirc)$ is the algebra of all continuous function on
the circle and if $\phi_1$ is given by integration with respect to Haar
measure, then $A$ is simple and purely infinite provided only that $A_2\ne\Cpx$
and $\phi_2$ is faithful but not a trace.

Parts of this work were done while the first named author took a part in the 
program "Quantenergodentheorie" organized by Walter Thirring at the Erwin 
Schr\"odinger International Institute for Mathematical Physics.
She would like to thank the Institute, the Institut f\"ur Theoretische 
Physik at Universit\"at Wien, and Professors Heide Narnhofer and Walter
Thirring for their warm hospitality.

\proclaim{Notation}\rm
We begin with some notation, which has appeared elsewhere.
Given an algebra $\Afr$ and subsets $S_\iota\subseteq \Afr$, ($\iota\in I$) let
$\Lambdao((S_\iota)_{\iota\in I})$ be the set of all words $w=a_1a_2\cdots a_n$
where $n\ge1$, $a_j\in S_{\iota_j}$ and
$\iota_1\ne\iota_2,\,\iota_2\ne\iota_3,\,\ldots,\,\iota_{n-1}\ne\iota_n$.
We will refer to the elements $a_1,\ldots,a_n$ as the letters of the word $w$;
we will sometimes regard the word as a product of specific letters, and
sometimes as an actual element of the algebra $\Afr$, as it suits the
situation.

Moreover, if a C$^*$--algebra $A$ and a state $\phi:A\to\Cpx$ are specified, we
will denote by $A\oup$ the kernel of $\phi$.
\endproclaim

\proclaim{Theorem}
Let $A_1$ be a reduced crossed product C$^*$--algebra,
$A_1=B\rtimes_{\alpha,r} G$,
where $G$ is an infinite discrete group and where $B$ is a unital
C$^*$--algebra.
Denote by $u_g$, ($g\in G$) the unitaries in $A_1$ arising from the reduced
crossed product construction and implementing the automorphisms $\alpha_g$ on
$B$.
Let $\phi_1$ be a faithful state on $B$ that is preserved by all the
automorphisms $\alpha_g$ and denote also by $\phi_1$ its extension to the state
on $A_1$ that vanishes on the subspace $Bu_g$ for every nontrivial $g\in G$.
Let $A_2$ be a unital C$^*$--algebra, $A_2\ne\Cpx$, with a faithful state
$\phi_2$;
let
$$ (A,\phi)=(A_1,\phi_1)*(A_2,\phi_2) $$
be the reduced free product of C$^*$--algebras.
Suppose that at least one of $\phi_1$ and $\phi_2$ is not a trace.

Then $A$ is purely infinite and simple.
\endproclaim
\demo{Proof}
Our strategy will be to show that $A$ is itself the reduced crossed product of
a C$^*$--subalgebra $D$ by the group $G$, where $D$ is (isomorphic to) the
reduced free product of infinitely many C$^*$--algebras;
a result from~\cite{\DykemaRordamZZProjII} will thereby show that $D$ is purely
infinite and simple.
We will then show that the action
of $G$ on $D$ is properly outer;
a result of Kishimoto and
Kumjian~\cite{\KishimotoKumjian} will thereby imply that $A$ is purely infinite
and simple.

\proclaim{Claim 1}
The family
$$ \bigl(B,(u_g^*A_2u_g)_{g\in G}\bigr) $$
is free with respect to $\phi$.
\endproclaim
\demo{Proof}
We must show that
$$ \Lambdao\bigl(\Bo,(u_g^*A_2\oup u_g)_{g\in G}\bigr)\subseteq\ker\phi.
\tag{\BugA} $$
Let $x$ be a word belonging to the left--hand--side of~(\BugA).
Splitting off the unitaries $u_g^*$ and $u_g$ from the letters in $x$,
then grouping together any neighbors in the resulting word belonging to $A_1$
and using that $u_{g_1}\Bo u_{g_2}^*\subseteq A_1\oup$ whenever $g_1,g_2\in G$
and that $u_{g_1}u_{g_2}^*\in A_1\oup$ if $g_1\ne g_2$, we see that $x$ is
equal to a word, $x'\in\Lambdao(A_1\oup,A_2\oup)$.
Hence $x\in\ker\phi$ by freeness.
This finishes the proof of Claim~1.
\enddemo

Let $D$ be the C$^*$--subalgebra of $A$ generated by
$B\cup\bigcup_{g\in G}u_g^*A_2u_g$.

\proclaim{Claim 2}
$D$ is simple and purely infinite.
\endproclaim
\demo{Proof}
Since $A_2\ne\Cpx$ there is a self--adjoint element, $x\in A_2\backslash\Cpx1$.
Let $\mu$ be the distribution of $x$;
namely, $\mu$ is the probability measure whose support is the spectrum of $x$
and such that $\phi_2(x^k)=\int_\Reals t^k\dif\mu(t)$ for all $k\ge1$.
A consequence of Bercovici and Voiculescu's
result~\scite{\BercoviciVoiculescuZZSuperconv}{Prop\. 8} is that for some $n$
large enough, the measure arising as the $n$--fold additive free convolution
$$ \mu_n\eqdef\undersetbrace{n\text{ times}}
\to{\mu\boxplus\mu\boxplus\cdots\boxplus\mu} $$
has support equal to an interval $[a,b]$ and is absolutely continuous with
respect to Lebesgue measure.
If $g_1,g_2,\ldots,g_n$ are distinct elements of $G$, then by Claim~1 the
distribution of $y\eqdef\sum_{j=1}^nu_{g_j}^*xu_{g_j}$ is $\mu_n$;
therefore $y$ generates an abelian subalgebra of
$$ D(g_1,\ldots,g_n)\eqdef C^*\Bigl(\bigcup_{j=1}^nu_{g_j}^*A_2u_{g_j}\Bigr) $$
on which $\phi$ is given by a measure without atoms;
it follows from~\scite{\DykemaHaagerupRordam}{Prop\. 4.1} that
$D(g_1,\ldots,g_n)$ contains a unitary $v$ satisfying $\phi(v)=0$;
(in fact, this proposition gives $\phi(v^k)=0$ for all nonzero integers $k$,
but we will not need this).
Therefore, partitioning the family $(u_g^*A_2 u_g)_{g\in G}$ into
subcollections of cardinality $n$, and including $B$ in one of these
subcollections, we see that $D$ is isomorphic to the free product of infinitely
many C$^*$--algebras with respect to faithful states,
$$ (D,\phi)\cong\freeprod{k=1}\infty(D_k,\psi_k), $$
where each $D_k$ contains a unitary that evaluates to zero under $\psi_k$.
Moreover, since either $\phi_2$ or $\phi_1\restrict_B$ is not a trace, at least
one of the $\psi_k$ is not a trace.
By~\scite{\DykemaRordamZZProjII}{Thm\. 2.1}, $D$ is therefore simple and purely
infinite.
This finishes the proof of Claim~2.
\enddemo

\proclaim{Claim 3}
$D$ has trivial relative commutant in $A$.
\endproclaim
\demo{Proof}
Let
$$ D_0=C^*\Bigl(\bigcup_{g\in G}u_g^*A_2u_g\Bigr)\subseteq D; $$
we will show that $D_0$ has trivial relative commutant in $A$, which will imply
the same for $D$.
Suppose that $x\in A$ and $x$ commutes with $D_0$;
our goal is to show that $x$ must belong to $\Cpx 1$.
Let $x_0=x-\phi(x)1$ and suppose, to obtain a contradiction, that $x_0\ne0$.
Since $\phi$ is faithful, $\nm{x_0}_2=\phi(x_0^*x_0)^{1/2}>0$.
Choose $\epsilon$ so that $0<\epsilon<\frac{\nm{x_0}_2}3$.
Since
$$ \Cpx1+\lspan\,
\Lambdao\Bigl(\Bo\cup\bigcup_{g\in G\backslash\{e\}}Bu_g,\;A_2\oup\Bigr)
\tag{\AdenseBug} $$
is a dense $*$--subalgebra of $A$, and since
$\Lambdao(\Bo\cup\bigcup_{g\in G\backslash\{e\}}Bu_g,A_2\oup)
\subseteq\ker\phi$,
there is a sum of finitely many words, $y=w_1+w_2+\cdots+w_m$ with
$w_1,w_2,\ldots,w_m\in
\Lambdao(\Bo\cup\bigcup_{g\in G\backslash\{e\}}Bu_g,A_2\oup)$,
such that
$\nm{x_0-y}<\epsilon$.
Let $F$ be the finite subset of $G$ whose elements are the identity element and
all nontrivial elements $g\in G$ for which some $w_j$ has a letter coming from
$Bu_g$.
From the proof of Claim~2, there is $n\in\Nats$ such that for any $n$ distinct
elements, $g_1,g_2,\ldots,g_n$ of $G$, there is a unitary
$$ v\in D(g_1,g_2,\ldots,g_n)=C^*\Bigl(\bigcup_{j=1}^nu_{g_j}^*A_2u_{g_j}\Bigr)
$$
with $\phi(v)=0$.
We take this unitary $v$ having ensured that the $n$ distinct elements satisfy
$g_j\notin F$ and $g_j^{-1}\notin F$ for every $j\in\{1,\ldots,n\}$.

Let us show that $vy$ and $yv$ are orthogonal with respect to the inner
product on $A$ induced by $\phi$, i.e\. that
$\langle yv,vy\rangle_\phi=\phi(v^*y^*vy)=0$.
Since
$$ \Cpx1+\lspan\,\Lambdao(u_{g_1}^*A_2\oup u_{g_1},u_{g_2}^*A_2\oup u_{g_2},
\ldots,u_{g_n}^*A_2\oup u_{g_n}) $$
is a dense $*$--subalgebra of $D(g_1,\ldots,g_n)$ and since (as can be seen
using Claim~1)
$$ \Lambdao(u_{g_1}^*A_2\oup u_{g_1},u_{g_2}^*A_2\oup u_{g_2},
\ldots,u_{g_n}^*A_2\oup u_{g_n})\subseteq\ker\phi, $$
for every $\eta>0$ there is a sum of finitely many words,
$z=w_1'+w_2'+\cdots+w_p'$ with
$$ w_1',\ldots,w_n'\in\Lambdao(u_{g_1}^*A_2\oup u_{g_1},
u_{g_2}^*A_2\oup u_{g_2},\ldots,u_{g_n}^*A_2\oup u_{g_n}), $$
such that $\nm{v-z}<\eta$.
But we see that each $w_j'$ is equal to a word
$$ w_j''\in\Lambdao\bigl(\bigl\{u_g\mid g\in G\backslash\{e\}\bigr\},
 A_2\oup\bigr) $$
where $w_j''$ begins with $u_{g_j^{-1}}$ and ends with $u_{g_k}$ some
$j,k\in\{1,\ldots,n\}$, and where $w_j''$ has length at least three.
Since
$$ w_1,\ldots,w_m\in\Lambdao\Bigl(\Bo\cup\bigcup_{g\in F\backslash\{e\}}Bu_g,
\;A_2\oup\Bigr), $$
when we consider a product $(w_{i_1}'')^*w_{j_1}^*w_{i_2}''w_{j_2}$ for
arbitrary $i_1,i_2\in\{1,\ldots,p\}$ and $j_1,j_2\in\{1,\ldots,m\}$, the choice
of the elements $g_1,\ldots,g_n$ ensures that there is
not too much cancellation and we are left with a reduced word
$$ (w_{i_1}'')^*w_{j_1}^*w_{i_2}''w_{j_2}=w\in
 \Lambdao\Bigl(\Bo\cup\bigcup_{g\in G\backslash\{e\}}Bu_g,A_2\oup\Bigr); $$
hence $\phi((w_{i_1}'')^*w_{j_1}^*w_{i_2}''w_{j_2})=0$.
This implies that $\phi(z^*y^*zy)=0$.
Since $\eta>0$ was arbitrary and
$|\phi(v^*y^*vy)-\phi(z^*y^*zy)|\le\eta(2+\eta)\nm y^2$,
we have $\phi(v^*y^*vy)=0$, i.e\. $yv$ and $vy$ are orthogonal.

We now obtain the contradiction.
Since $x_0$ belongs to the commutant of $D_0$, we must have $vx_0-x_0v=0$.
But by orthogonality of $vy$ and $yv$,
$$ \nm{vy-yv}\ge\nm{vy-yv}_2>\nm{vy}_2=\nm y_2 $$
and hence
$$ \nm{vx_0-x_0v}\ge\nm{vy-yv}-2\epsilon>\nm y_2-2\epsilon\ge
\nm{x_0}_2-3\epsilon>0, $$
which is a contradiction.
This finishes the proof of Claim~3.
\enddemo

\proclaim{Claim 4}
For every nontrivial $g\in G$, $\beta_g\eqdef\Ad(u_g)$ is an outer
automorphism of $D$, $g\mapsto\beta_g$ is a group homomorphism and $A$ is
isomorphic to the reduced crossed product $D\rtimes_{\beta,r}G$.
\endproclaim
\demo{Proof}
Clearly, $\Ad(u_g)$ is an automorphism of $D$, for every $g\in G$ and
$g\mapsto\beta_g$ is a group homomorphism.
From the density of~(\AdenseBug) in $A$ and the fact that $u_gB=Bu_g$, we see
that $\lspan\bigcup_{g\in G}Du_g$ is dense in $A$.
Moreover, whenever $g'\in G$ is nontrivial, $Du_{g'}\subseteq\ker\phi$;
this can be seen by approximating an arbitrary element of $Du_{g'}$ by sums of
words each belonging to
$\{u_{g'}\}\cup\Lambdao(\Bo,(u_g^*A_2\oup u_g)_{g\in G})u_{g'}$.
As the GNS representation of $\phi$ is faithful on $A$, one sees that $A$ is
isomorphic to the reduced crossed product $D\rtimes_{\beta,r}G$.

We will now show that $\beta_g$ is an outer automorphism of $D$ whenever
$g\ne e$.
Indeed, if it were inner then letting $v_g\in D$ be such that
$\beta_g=\Ad(v_g)$, we would have $u_g^*v_g$ commuting with $D$.
By Claim~3, this would imply that $u_g$ is a scalar multiple of $v_g$, hence
belongs to $D$,
which contradict that $Du_g\subseteq\ker\phi$.
This finishes the proof of Claim~4.
\enddemo

Now that $A$ is seen to be the crossed product of a simple, purely infinite
C$^*$--algebra by an infinite discrete group acting by outer automorphisms,
Kishimoto and Kumjian's result~\scite{\KishimotoKumjian}{Lemma 10} shows that
$A$ is simple and purely infinite.
\QED


\Refs

  \ref \no \Avitzour \by D\. Avitzour \paper Free products of C$^*$--algebras
    \jour Trans\. Amer\. Math\. Soc\. \vol 271 \yr 1982 \pages 423-465 \endref

  \ref \no \BercoviciVoiculescuZZSuperconv \by H\. Bercovici, D\. Voiculescu
    \paper Superconvergence to the free central limit theorem and failure of
    Cram\'er theorem for free random variables
    \jour Prob\. Theory Relat\. Fields
    \vol 102 \yr 1995 \pages 215-222 \endref

  \ref \no \CuntzZZKthCertain \by J\. Cuntz
    \paper K--theory for certain C$^*$--algebras
    \jour Ann\. of Math\. \vol 113 \yr 1981 \pages 181-197 \endref

  \ref \no \DykemaZZFaithful \manyby K.J\. Dykema
    \paper Faithfulness of free product states
    \jour J\. Funct\. Anal\. \vol 154 \yr 1998 \pages 223-229 \endref

  \ref \no \DykemaZZPII \bysame
    \paper  Purely infinite simple $C^*$-algebras arising from free product
    constructions, II
    \jour Math\. Scand\. \toappear \endref

  \ref \no \DykemaHaagerupRordam \by K.J\. Dykema, U\. Haagerup, M\. R\o rdam
    \paper The stable rank of some free product C$^*$--algebras
    \jour Duke Math\. J\. \vol 90 \yr 1997 \pages 95-121
  \moreref \paper correction \vol 94 \yr 1998 \page 213 \endref

  \ref \no \DykemaRordamZZPI \manyby K.J\. Dykema, M\. R\o rdam
    \paper  Purely infinite simple $C^*$-algebras arising from free product
    constructions
    \jour Can\. J\. Math\. \vol 50 \yr 1998 \pages 323-341 \endref

  \ref \no \DykemaRordamZZProj \bysame
    \paper Projections in free product C$^*$--algebras
    \jour Geom\. Funct\. Anal\. \vol 8 \yr 1998 \pages 1-16 \endref

  \ref \no \DykemaRordamZZProjII \bysame
    \paper Projections in free product C$^*$--algebras, II
    \jour Math\. Z\. \toappear \endref

  \ref \no \KishimotoKumjian \by A\. Kishimoto, A\. Kumjian
    \paper Crossed products of Cuntz algebras by quasi--free automorphisms
    \jour Fields Inst\. Commun\. \vol 13 \yr 1997 \pages 173-192 \endref

  \ref \no \VoiculescuZZSymmetries \by D\. Voiculescu
    \paper Symmetries of some
    reduced free product C$^{\ast}$--algebras
    \inbook Operator Algebras and Their Connections with Topology and Ergodic
    Theory
    \bookinfo Lecture Notes in Mathematics \vol1132
    \publ Springer--Verlag \yr 1985
    \pages 556--588 \endref

  \ref \no \VDNbook \by D\. Voiculescu, K.J\. Dykema, A\. Nica
    \book Free Random Variables \bookinfo CRM Monograph Series vol\.~1
    \publ American Mathematical Society \yr 1992 \endref

\endRefs
\enddocument